\documentclass[12pt]{article}
\usepackage{amsfonts}
\usepackage[dvips]{graphics}
\newtheorem{theorem}{Theorem}

\newtheorem{definition}{Definition}

\begin{document}
\title {A Liouville comparison principle for entire sub- and super-solutions of
the equation $u_t-\Delta_p (u) = |u|^{q-1}u$.}
\author{Vasilii V. Kurta.}
\maketitle
\begin{abstract}
\noindent We establish a Liouville comparison principle for entire
sub- and super-solutions of the equation $(\ast)$ $w_t-\Delta_p (w)
= |w|^{q-1}w$ in the half-space  ${\mathbb S}= {\mathbb R}^1_+\times
{\mathbb R}^n$, where $n\geq 1$, $q>0$ and  $ \Delta_p
(w):=\mbox{div}_x\left(|\nabla_x w|^{p-2}\nabla_x w\right)$,
$1<p\leq 2$. In our study we impose neither  restrictions on the
behaviour of entire  sub- and super-solutions on the hyper-plane
$t=0$,  nor any  growth conditions on the behavior   of them or any
of their partial derivatives at infinity. We prove that if $1<q\leq
p-1+\frac pn$, and $u$ and $v$ are, respectively, an entire weak
super- and an entire  weak sub-solution of ($\ast$) in $\Bbb S$
which belong, only locally in $\Bbb S$, to the corresponding Sobolev
space and are such that $u\leq v$, then $u\equiv v$. The result is
sharp. As direct corollaries we obtain both new and known
Fujita-type and Liouville-type results.
\end{abstract}

\thispagestyle{empty}
\section{Introduction and Definitions.}

The purpose of this work is to obtain a  Liouville comparison
principle of  elliptic type for entire   weak sub- and
super-solutions of the equation
\begin{eqnarray}
w_t-\Delta_p (w) = |w|^{q-1}w
\end{eqnarray}
in the half-space  ${\mathbb S}= (0,+\infty)\times {\mathbb R}^n$,
where $n\geq 1$ is a natural number, $q>0$ is a  real number and $
\Delta_p (w):=\mbox{div}_x\left(|\nabla_x w|^{p-2}\nabla_x
w\right)$, $1<p\leq 2$,  defines  the well-known $p$-Laplacian
operator. Under entire sub- and super-solutions of (1)  we
understand sub- and super-solutions of (1) defined in the whole
half-space $\mathbb S$ and under Liouville theorems of elliptic type
we understand Liouville-type theorems which, in their formulations,
have no restrictions on the behaviour of sub- or super-solutions to
the parabolic equation (1) on the hyper-plane $t=0$. We would also
like to underline that we impose no growth conditions on the
behavior of sub- or super-solutions of (1), as well as of all their
partial derivatives,  at infinity.

\begin{definition}
Let $n\geq 1$, $p>1$ and $q>0$. A  function $u=u(t,x)$ defined and
measurable in $\mathbb S$ is called  an entire   weak super-solution
of the equation (1) in $\mathbb S$ if it belongs to the function
space $ L_{q,loc}(\mathbb S)$, with $u_t \in L_{1, loc}(\mathbb S)$
and $|\nabla_x u|^{p}\in L_{1, loc}(\mathbb S)$, and satisfies the
integral inequality
\begin{eqnarray}
\int\limits_{\mathbb S}\left[u_t\varphi+\sum_{i=1}^n
|\nabla_x u|^{p-2}u_{x_i}\varphi_{x_i} -|u|^{q-1}u \varphi \right]dtdx
\geq 0
\end{eqnarray}
for every  non-negative function $\varphi \in  C^\infty (\mathbb S)$
with compact support in $\mathbb S$,  where $C^{\infty}({\mathbb
S})$ is the space of all functions defined and infinitely
differentiable in  ${\mathbb S}$.
\end{definition}
\begin{definition} A function $v=v(t,x)$ is an entire   weak sub-solution of (1)
if $u=-v$ is an entire   weak  super-solution of (1).
\end{definition}

\section{Results.}

\begin{theorem} Let $n\geq 1$,  $2\geq p> 1$
and  $1<q\leq p-1 +\frac {p}n$,  and let  $u$ be an entire   weak
super-solution and $v$ an entire   weak sub-solution of  (1) in
$\mathbb S$ such that $u\geq v$. Then $u \equiv v$ in $ \mathbb S$.
\end{theorem}

The result in Theorem 1, which evidently has a comparison principle
character, we term a Liouville-type comparison principle, since, in
the particular cases when $u\equiv 0$ or $v\equiv 0$, it becomes a
Liouville-type theorem of elliptic type, respectively,  for entire
sub- or super-solutions of (1).

Since in Theorem 1 we impose no conditions on the behaviour of
entire  sub- or super-solutions of  the equation (1) on the
hyper-plane $t=0$, we can formulate,  as a direct corollary  of the
result in Theorem 1, the following comparison principle, which in
turn one can term a Fujita comparison principle,  for entire sub-
and super-solutions of the Cauchy problem for the equation (1). It
is clear that in the particular cases when $u\equiv 0$ or $v\equiv
0$, it becomes a Fujita-type theorem, respectively, for entire sub-
or super-solutions of the Cauchy problem for the equation (1).

\begin{theorem} Let $n\geq 1$,  $2\geq p> 1$
and  $1<q\leq p-1 +\frac {p}n$, and let  $u$ be an entire   weak
super-solution and $v$ an entire   weak sub-solution of the Cauchy
problem, with possibly different initial data for $u$ and $v$, for
the equation (1) in $\mathbb S$ such that $u\geq v$. Then $u \equiv
v$ in $\mathbb S$.
\end{theorem}


Note that the results  in Theorems 1 and 2 are sharp. The sharpness
of these for $q>  p-1 +\frac {p}n\geq 1$ follows, for example, from
the existence of non-negative self-similar entire solutions to (1)
in $\Bbb S$, that was shown  in [1]. Also, there one can find a
Fujita-type theorem on blow-up of non-negative entire solutions of
the Cauchy problem for (1), which was obtained as a very interesting
generalization of the famous blow-up result in [2] to quasilinear
parabolic equations. For $0<q\leq1$, it is evident that the function
$u(t,x)=e^t$  is a positive entire classical super-solution of (1)
in  $\Bbb S$.

To prove the results in Theorems 1 and 2 we further develop the
approach that was proposed for solving similar problems for
semilinear parabolic equations in [3]. A new key point in our proof
is using the fact that for $1<p\leq 2$ the $p$-Laplacian operator
$\Delta_p$ satisfies the $\alpha$-monotonicity condition (see, e.g.,
[4]) with $\alpha=p$.

\vskip 10pt \noindent {\bf Acknowledgments.}

\vskip 10pt \noindent This research is financially supported by the
Alexander von Humboldt Foundation (AvH). The author is very grateful
to AvH  for the  opportunity to visit the Mathematical Institute of
K\"oln University and to Professor B. Kawohl for his cordial
hospitality during this visit.

\vskip 10pt

\noindent \textbf{Author's address:}

\vspace{10 mm}

\noindent Vasilii V. Kurta

\noindent Mathematical Reviews

\noindent 416 Fourth Street, P.O. Box 8604

\noindent Ann Arbor, Michigan 48107-8604, USA

\noindent \textbf {e-mail:} vkurta@umich.edu, vvk@ams.org

\vskip 10pt




\begin{thebibliography}{5}



\bibitem{} V.A. Galaktionov and H.A. Levine, \emph{ A general approach to
critical Fujita exponents in nonlinear parabolic problems},
Nonlinear Anal.  (1998) v. 34,  no. 7, 1005--1027.

\bibitem{}
H. Fujita, \emph{On the blowing up of solutions of the Cauchy
problem for $u_t=\Delta u + u^{1+\alpha}$}, J. Fac. Sci. Univ.
Tokyo, Sect. I (1966), v. 13, 109--124.


\bibitem{} A.G. Kartsatos and V.V. Kurta, \emph{ On a comparison
principle and the critical Fujita exponents for solutions of
semilinear parabolic inequalities}, J. London Math. Soc. (2) 66
(2002), no. 2, 351--360.



\bibitem{}
V.V. Kurta, \emph{Comparison principle for solutions of parabolic
inequalities}, C. R. Acad. Sci. Paris, S\'{e}rie I, 322 (1996),
1175--1180.





\end{thebibliography}
\end{document}